\documentclass[11pt]{amsart}

\usepackage{latexsym, amssymb,lscape,epsfig}

\theoremstyle{plain}
\newtheorem{theorem}{Theorem}[section]
\newtheorem{proposition}{Proposition}[section]
\newtheorem{conjecture}{Conjecture}[section]

\theoremstyle{definition}
\newtheorem{lem}{Lemma}[section]

\newtheorem{example}{Example}[section]
\newtheorem{exercise}{Exercise}

\newenvironment{renumerate}%
{%
\begin{enumerate}}%
{\end{enumerate}%
}%

\newenvironment{theo}[1]%
{\begin{theorem}\label{T:#1}}%
{\end{theorem}}

\newenvironment{prop}[1]%
{\begin{proposition}\label{T:#1}}%
{\end{proposition}}

{\begin{conjecture}\label{T:#1}}%
{\end{conjecture}}

\newenvironment{definition}%
{\vskip6pt%
\noindent%
{\bf Definition.}}%
{\vskip6pt}

\newenvironment{remark}%
{\vskip6pt%
\noindent%
{\it Remark.}}%
{\vskip6pt}

{\vskip6pt%
\noindent%
{\it Remarks}. %
\begin{renumerate}}%
{\end{renumerate}\vskip6pt}

{\begin{exercise}\label{T:#1}}%
{\end{exercise}}

\newenvironment{ex}[1]%
{\begin{example}\label{T:#1}}%
{\end{example}}

\newcommand{\C}{\text{$\mathbb C$}}
\newcommand{\N}{\text{$\mathbb N$}}

\renewcommand{\frak}[1]{\text{$\mathfrak{#1}$}}
\newcommand{\J}{\text{$\mathcal{J}$}}

\newcommand{\ga}{\text{$\alpha$}}

\newcommand{\e}{\text{$\varepsilon$}}

\newcommand{\gO}{\text{$\Omega$}}
\newcommand{\go}{\text{$\omega$}}
\newcommand{\gf}{\text{$\varphi$}}

\newcommand{\del}{\text{$\partial$}}

\newcommand{\delbar}{\text{$\overline{\partial}$}}

\newcommand{\im}{\mathrm{Im}\,}

\newcommand{\mc}[1]{\text{$\mathcal{#1}$}}

\newcommand{\into}{\rightarrow}
\newcommand{\noqed}{\let\qed\relax}

\newcommand{\Gg}{\mathfrak{g}}

\numberwithin{equation}{section}

\begin{document}
\pagestyle{headings}

\title{The Lefschetz property, formality and blowing up in symplectic geometry}
\author{Gil Ramos Cavalcanti}
\address{Mathematical Institute, St. Giles 24-29, Oxford, OX1 3BN, UK}
\email{gil.cavalcanti@new.oxon.org}
\thanks{Researcher supported by CAPES (Coordena\c c\~ao de Aperfei\c coamento de Pessoal de N\'ivel Superior, Minist\'erio da Educa\c c\~ao e Cultura), Brazilian Government. Grant 1326/99-6}

\begin{abstract}
In this paper we study the behaviour of the Lefschetz property under the blow-up construction. We show that it is possible to reduce the dimension of the kernel of the Lefschetz map if we blow up along a suitable submanifold satisfying the Lefschetz property. We use that, together with results about Massey products, to construct compact nonformal symplectic manifolds satisfying the Lefschetz property.
\end{abstract}

\maketitle
\markboth{\sc G. R. Cavalcanti}{\sc The Lefschetz property, formality and blowing up}

\section*{Introduction}

In \cite{Kos85} Koszul introduced an operator $\delta$ for Poisson manifolds which consists of the exterior derivative twisted by the Poisson bivector. This operator was further studied by Brylinski \cite{Br88} who pointed some similarities between $\delta$ and $d^*$, the Riemannian adjoint of $d$. He called forms that are both $d$ and $\delta$ closed {\it symplectic harmonic}. Later, Yan \cite{Ya96} and Mathieu \cite{Ma90} proved independently that the existence of a harmonic representative in each cohomology class is equivalent to the strong Lefschetz property (or just Lefschetz property, for short).

\vskip6pt
\noindent
{\bf Lefschetz property.} A symplectic manifold $(M^{2n},\go)$ satisfies the {\it Lefschetz property at level $k$} if the map
$$[\go^{n-k}] : H^{k}(M) \into H^{2n-k}(M)$$
is surjective. It satisfies the {\it Lefschetz property} if these maps are surjective for $0 \leq k \leq n$.
\vskip6pt

Recently, Merkulov \cite{Me98} proved that, in a compact symplectic manifold, the existence of symplectic harmonic forms on each cohomology class (and therefore the Lefschetz property) is equivalent to the $d\delta$-lemma:

\vskip6pt
\noindent
{\bf $d\delta$-lemma.} A symplectic manifold {\it satisfies the $d\delta$-lemma} if
$$ \im d \cap \ker \delta =\im \delta \cap \ker d = \im d\delta.$$
\vskip6pt

\newpage

Therefore we have the following implications for a compact symplectic manifold:
\vskip6pt
\setlength{\unitlength}{1mm}
\begin{picture}(60,10)(-65,-5)
\linethickness{1pt}
\thinlines
\put(-50,0){\makebox(0,0){Lefschetz property}}
\put(5,0){\makebox(0,0){Harmonic representatives}}
\put(55,0){\makebox(0,0){$d\delta$-lemma}}
\put(35,0){\makebox(0,0){$\Leftrightarrow$}}
\put(-25,0){\makebox(0,0){$\Leftrightarrow$}}
\end{picture}

Although originally thought of as an analogue of $d^*$, the existence of a $d\delta$-lemma draws similarities between $\delta$ and the complex $d^c = -i(\del - \delbar)$. More evidence of the connection between these operators comes from generalized complex geometry, a theory recently introduced by Hitchin \cite{Hi03} and further developed by Gualtieri \cite{Gu03}. As pointed out by Gualtieri, both $d^c$ and $\delta$ are particular cases of a general operator $d^{\J}$ on a generalized complex manifold and implications of the $dd^{\J}$-lemma have been studied in the author's thesis \cite{Ca04b}.

In this paper we investigate whether two particular properties of the standard $dd^c$-lemma can be translated to the symplectic case. Firstly, it is a result of Parshin \cite{Pa66} that the $dd^c$-lemma is preserved by holomorphic rational equivalences, in particular, by the blow-up along a complex submanifold. Secondly, it implies that the manifold is formal in the sense of Sullivan \cite{DGMS75}.

As far as the first point is concerned, we study how the Lefschetz property, and hence the $d\delta$-lemma, behaves under the operation of symplectic blow-up introduced by McDuff \cite{McD84}. McDuff used symplectic blow-up to give the first simply-connected example of a non-K\"ahler symplectic manifold, the blow-up of $\C P^5$ along a symplectically embedded Thurston manifold (the nilmanifold with structure $(0,0,0,12)$). This example of McDuff fails to be K\"ahler, since it does not satisfy the Lefschetz property. This means that even if the ambient space satisfies the $d\delta$-lemma, the blown-up manifold may not do the same. 

In this paper we study systematically how the Lefschetz property behaves under the blow-up, in particular we seek conditions under which we can assure that the blown-up manifold will satisfy the Lefschetz property. We prove that this is the case if both submanifold and ambient manifold satisfy the property and the co-dimension is high enough (see Theorem \ref{T:lefschetz}). Moreover we study the blow-down map and show that even if the blown-up manifold satisfies the Lefschetz property, the original ambient manifold will not necessarily do so (see Theorem \ref{T:lefschetz in m2} and Proposition \ref{T:i=2d}). Together with McDuff's example, this shows that we can not decide whether the blow-up will or will not satisfy the Lefschetz property based solely on the ambient manifold. 

The second purpose of this paper is to answer the question of correlation of the $d\delta$-lemma and formality. One of the implications is known to be false: Gompf produced a simply-connected 6-manifold which does not satisfy the Lefschetz property \cite{Go95}, but which is formal by Miller's result \cite{Mi79}. The converse implication had been conjectured by Babenko and Taimanov \cite{BT00a} and was the object of study of other papers \cite{IRTU03,LO94}.

Our starting point in this case is that Babenko and Taimanov studied thoroughly the behaviour of Massey products under blow-up \cite{BT00a,BT00b} and these tend to `survive' in the blow-up, which is markedly different from the behaviour of the Lefschetz property. Using this approach, we produce an example of a compact nonformal symplectic manifold satisfying the Lefschetz property by blowing-up a 6-nilmanifold along a suitable torus. We also produce a 4-dimensional example using Donaldson submanifolds \cite{Do96} and results of Fern\'andez and Mu\~noz \cite{FM02} and, with a further blow-up, we obtain a simply-connected 12-dimensional nonformal compact symplectic manifold satisfying the Lefschetz property.

This paper is organized as follows. In the first section we explain briefly how the blow-up is done in symplectic geometry and derive the cohomology algebra of the blow-up from the cohomologies of the ambient manifold and submanifold as well as the Chern and Thom classes of the normal bundle of the submanifold. In Section \ref{S:lefschetz}, we study how the Lefschetz property behaves under blow-up, initially in the case of the blow-up along an embedded surface and later the general case. In Section \ref{S:massey} we recall the definition of Massey products, their relation with formality and their behaviour under blow-up. In Section \ref{S:nonformal examples} we provide examples of compact nonformal symplectic manifolds satisfying the Lefschetz property.

The author has been informed that other examples of compact nonformal symplectic manifolds satisfying the Lefschetz property have been found independently by Amor\'os and Kotschick using different methods.

{\bf Acknowledgments.} I thank Oliver Thomas for comments on determinants which were useful in lemma \ref{T:oliver} and Vicente Mu\~noz for the argument used in Example \ref{T:munoz}. I also thank Marisa Fern\'andez and Marco Gualtieri for helpful suggestions and  Nigel Hitchin for inspiring discussions, guidance and help with the editing of the text.

\section{The Symplectic Blow-up}\label{S:blowup}

We begin by giving a description of the cohomology ring of the blown-up
manifold in terms of the cohomology rings of the ambient manifold and the
embedded submanifold and the Chern and Thom classes of the normal
bundle of the embedding. We shall outline the blow-up construction in
order to fix some notation. For a detailed presentation we refer to
\cite{McD84}.

Assume that $i:(M^{2d},\sigma) \hookrightarrow
(X^{2n},\go)$ is a symplectic
embedding, with $M$ compact. Let $k=n-d$. In these circumstances we can 
choose a
complex structure in $TX$ that restricts to one in $TM$, and hence
also to the normal bundle $E
\stackrel{\pi}{\rightarrow} M$. Therefore $E$ is a complex
bundle over $M$ and one can form its projectivization
$$ \mathbb{C}P^{k-1} \longrightarrow \tilde{M} \longrightarrow M$$
and also form the ``tautological'' line bundle $\tilde{E}$ over $\tilde{M}$:
 the subbundle of $\tilde{M}\times E$ whose fibers are the
elements $\{([v],\lambda v), \lambda \in \mathbb{C}\}$. We have
the following commutative diagram
\begin{equation}\label{E:diagram}
\begin{array}{ccccc}
\tilde{E_0} & \longrightarrow &\tilde{E} &
\stackrel{q}{\longrightarrow}& \tilde{M}\\
\Big\downarrow \vcenter{\rlap{$\gf$}}& & \Big\downarrow \vcenter{\rlap{$\gf$}}& &
\Big\downarrow \vcenter{\rlap{$p$}}\\
E_0 & \longrightarrow &(E,\go) &
\stackrel{\pi}{\longrightarrow}& (M,\sigma)
\end{array}
\end{equation}
where $q$ and \gf\ are the projections over $\tilde{M}$ and $E$
respectively, $E_0$ is the complement of the zero section in $E$ and
$\tilde{E}_0$ the complement of the zero section in $\tilde{E}$.

It is easily seen, $E_0$
and $\tilde{E}_0$ are diffeomorphic via \gf. Furthermore, if we let $V$
be a sufficiently small disc subbundle in $E$ with its canonical
symplectic structure \go, then it is symplectomorphic to a neighbourhood
of $M \subset X$ and we identify the two from now on. Letting
$\tilde{V} = \gf^{-1}(V)$,
we can form the manifold
$$ \tilde{X} = \overline{X-V} \cup_{\del V}\tilde{V}.$$
Then, the map \gf\ can be extended to a map  $f:\tilde{X}
\into X$, being the identity in the complement of $\tilde{V}$. The manifold $\tilde{X}$ is the blow-up of $X$ along $M$ and $f:\tilde{X}
\into X$ is the projection of the blow-up or the blow-down map.

\begin{lem} (McDuff \cite{McD84})
There is a unique class $a \in H^2(\tilde{M})$ which restricts to the standard K\"ahler class on each fiber of $\tilde{M} \into M$ and pulls back to the trivial class in $\tilde{E}_0$. Moreover, $H^{\bullet}(\tilde{E}) \cong H^{\bullet}(\tilde{M})$ is a free module over $H^{\bullet}(M)$ with generators $1,a,\cdots,a^{k-1}$. 
\end{lem}

\begin{theo}{cohomology}
{\em (McDuff \cite{McD84})} If the codimension of $M$ is at least 4, the
fundamental groups of $X$ 
and the blown up manifold $\tilde{X}$ are isomorphic. Further, there
is a short exact sequence
$$0 \into H^*(X) \rightarrow H^*(\tilde{X}) \into A^* \into 0,$$
where $A^*$ is the free module over $H^*(M)$ with generators $a,
\cdots, a^{k-1}$.
Moreover, there is a representative $\ga$ of $a$ with support in the 
tubular neighbourhood $V$ such that, for \e\ small enough, the form 
$\tilde{\go} = f^*(\go) + \e \ga$ is a symplectic form in $\tilde{X}$.
\end{theo}

\begin{remark}
The symplectic structure in the blow-up, $\tilde{X},$ is {\it not} determined by the one in $X$. The presence of the parameter $\e$ in the symplectic form is just one indicative that this is the case. Futhermore, for each $\e$, the symplectic structure in $\tilde{X}$ also depends on other choices made during the construction, as the almost complex structure taming $\go$ and the identification of the normal bundle with the tubular neighbouhood. See \cite{MS00}, page 231.
\end{remark}

\begin{remark}
If the submanifold has codimension less than 4, then the blow-up will be just $X$ again and hence the theorem is trivially true.
\end{remark}

\begin{remark}
As is observed by McDuff \cite{McD84}, the Leray-Hirsch theorem implies that $a^k$ is related to $a, \cdots, a^{k-1}$ in $\tilde{E}$ by
$$
a^k = -c_k -c_{k-1} a + \cdots - c_1 a^{k-1},
$$
where the $c_j$'s are the Chern classes of the normal bundle $E$.

In \cite{RT00} it is shown that in $\tilde{X}$ this relation becomes
$$
a^k = -f^*(t) -c_{k-1} a + \cdots - c_1 a^{k-1},
$$
where $t$ is the Thom class of the embedding $M \hookrightarrow X$ and $f:\tilde{X} \into X$ the projection of the blow-up.
\end{remark}

With this, we have a complete description of the cohomology ring of
$\tilde{X}$. For $v_1, v_2 \in H^*(X)$ and $u_1, u_2 \in H^*(M)$,
\begin{equation}\label{E:product rules}
\begin{cases}
f^*(v) \wedge f^*(w) &= f^*(v\wedge w);\\
f^*(v) a &= i^*(v) a;\\
a u_1 \wedge a u_2 &= a^2 u_1 \wedge u_2\\
a^k &= -f^*(t) -c_{k-1} a + \cdots - c_1 a^{k-1};\\
f^*(t) \wedge u_1 &=  f^*(t \wedge u_1), \mbox{ the Thom map extended 
to } X.
\end{cases}
\end{equation}

\section{The Lefschetz Property and Blowing-up}\label{S:lefschetz}

Now we move on to study how the Lefschetz property behaves under blow-up. The first case to look at would be the blow-up of a point, but, as
we will see, this does not change the kernel of the Lefschetz map at
any level (cf. Theorem \ref{T:i>2d}). The next case would be a
surface. Here, on the one hand, the situation is simple enough for us
to be able to give a fairly complete account of what happens, and, on the other, we can already see that in this case it is possible to decrease the
dimension of the kernel of the Lefschetz map.

\subsection{Blowing up along a Surface}
Assume that $i:(M^2,\sigma) 
\hookrightarrow
(X^{2n},\go)$ is 
a surface symplectically embedded in $X$, $M$ and $X$ are compact, and
let $\tilde{X}$ be the 
blow-up of $X$ along $M$. In $H^1(\tilde{X})$ things go as follows
\begin{equation}\label{E:h1}
\begin{aligned}
(f^*(\go) + \e a)^{n-1} f^*(v) &= f^*(\go^{n-1}v) + \e^{n-1}
a^{n-1}i^*v \\
&= f^*(\go^{n-1}v - \e^{n-1} t v)
\end{aligned}
\end{equation}
and if Lefschetz holds for $X$ and $\e$ is small enough Lefschetz will
also hold for $\tilde{X}$, or, more generally,
$\dim(\ker(\tilde{\go}^{n-1})) \leq  
\dim(\ker(\go^{n-1}))$. Now we proceed to show that in certain 
conditions the inequality holds.

\begin{lem}\label{T:lemma1}
Let $i: (M^2,\sigma) \hookrightarrow (X^{2n},\go)$ be a symplectic
embedding, $M$ and $X$ be compact and $t$ be the Thom class of this embedding. The following
are equivalent: 
\begin{enumerate}
\item There are $v_1, v_2 \in H^1(X)$ in $\ker(\go^{n-1})$
such that $i^*(v_1\wedge v_2) \neq 0$;
\item There exists $v_1 \in \ker(\go^{n-1})$ such that $t\wedge v_1 \not
\in \im(\go^{n-1})$.
\end{enumerate}
\end{lem}
\begin{proof}
Assuming (1), by the defining property of the Thom class,
$$ \int_X t \wedge v_1 \wedge v_2 = \int_M i^*(v_1 \wedge v_2) \neq 0,$$
and, since both $v_1$ and $v_2$ pair trivially with $\im(\go^{n-1})$,
but pair nontrivially with $t \wedge v_i$, we
see that $t \wedge v_i \not \in \im(\go^{n-1})$. So (1) implies (2).

On the other hand, assume that there is a $v_1$ satisfying (2). Let $\{a_i\}$ be a
basis for $\ker (\go^{n-1})$ and $\{\tilde{a}_i\}$ be a basis for a
complement. Since $(H^1)^* \cong H^{2n-1}$, we can view the dual basis
$\{a_i^*,\tilde{a}_i^*\}$ as a basis for $H^{2n-1}$. Then we note that
$\im(\go^{n-1}) \subset \mbox{span}\{\tilde{a}_i^*\}$, and since these spaces
have the same dimension they are the same. Therefore, the condition
$t\wedge v_1 
\not \in \im(\go^{n-1})$ implies that it pairs nontrivially with some
of the $a_i$. Let  $v_2$ be such an $a_i$. Then again by the defining
property of the Thom class we have
$$\int_M i^*(v_1 \wedge v_2) = \int_X t \wedge v_1 \wedge v_2 \neq 0,$$
and $i^*(v_1\wedge v_2) \neq 0$.
\end{proof}

\begin{lem}\label{T:h1}
If the equivalent conditions (1) and (2) of the previous lemma are
satisfied and \e\ is 
small enough, then
$$ \dim(\ker(\tilde{\go}^{n-1}:H^1(\tilde{X}) \into
H^{2n-1}(\tilde{X}))) \leq \dim(\ker(\go^{n-1}:H^1(X) \into
H^{2n-1}(X))) - 2.$$
\end{lem}
\begin{proof}
Let $V$ be a complement of $\ker(\go^{n-1})$ in $H^1(X)$ and
$v_1$ and 
$v_2$ the cohomology classes satisfying condition (1) of 
Lemma \ref{T:lemma1}. Then, since neither $t \wedge v_1$ or $t \wedge
v_2$ is in 
$\im(\go^{n-1})$, for \e\ small enough, equation \eqref{E:h1}
shows that
$f^*(t \wedge v_i) \not \in \im(\tilde{\go}^{n-1}|_V)$, 
since
$\tilde{\go}^{n-1}|_V$ is simply a perturbation of the injection 
$\go^{n-1}|_V$. On the other hand, 
$\tilde{\go}^{n-1} f^*(v_i) = -\e^{n-1} f^*(t\wedge v_i)$, and therefore $f^*(t
\wedge v_i)$ is in the image of $\tilde{\go}^{n-1}$, so
$$\dim(\im(\tilde{\go}^{n-1})) \geq \dim(\im(\go^{n-1})) + 2$$ 
and the result follows.
\end{proof}

Now we move on to $H^2(\tilde{X})$, where we have
\begin{equation}\label{E:h2}
\begin{aligned}
(f^*(\go) + \e a)^{n-2} (f^*(v_2) + a v_0) =&\\
= f^*(\go^{n-2}v_2  -  \e^{n-2} t v_0)& + \e^{n-3}a^{n-2}((n-2) \sigma v_0
 + \e( i^*v_2 - c_1 v_0))
\end{aligned}
\end{equation}
and then we observe that the map above is a perturbation of
$$f^*(v_2) + a v_0 \mapsto f^*(\go^{n-2}v_2) + \e^{n-3}a^{n-2}(n-2) \sigma v_0.$$
Therefore for \e\ small enough, Lefschetz will hold for $\tilde{X}$
if it holds for $X$, or more generally $\dim(\ker(\tilde{\go}^{n-2}))
 \leq \dim(\ker(\go^{n-2}))$.

Again, we may have the inequality.
\begin{lem}\label{T:lemma2}
Let $i: (M^2,\sigma) \hookrightarrow (X^{2n}\go)$ be a symplectic
embedding, $M$ and $X$ be compact and $t$ be the Thom class of this embedding. The following
are equivalent: 
\begin{enumerate}
\item There exists $v \in \ker \go^{n-2}$ such that $i^*v \neq 0$;
\item The Thom class $t$ is not in the image of $\go^{n-2}$.
\end{enumerate}
\end{lem}
\begin{proof}
The proof is the same as the one for Lemma
\ref{T:lemma1}. Assuming (1), by definition of the Thom class,
$$ \int_X t v = \int_M i^*v \neq 0.$$
On the other hand any vector in the kernel of $\go^{n-2}$ pairs
trivially with $\im(\go^{n-2})$, so $t\not \in \im(\go^{n-2})$.

Conversely, we let again $\{a_i\}$ be a basis for $\ker(\go^{n-2})$,
$\{\tilde{a}_i\}$ a basis for a complement and
$\{a^*_i,\tilde{a}^*_i\}$ the dual basis and again identify the
dual space with $H^{2n-2}$. Then we see that $\im(\go^{n-2}) =
\mbox{span}\{\tilde{a}^*_i\}$ and, since $t \not \in \im(\go^{n-2})$,
$t$ must pair nontrivially with at least one of the $a_i$'s. Call it $v$.
\end{proof}

\begin{lem}\label{T:h2}
If the equivalent conditions (1) and (2) of the previous lemma are
satisfied and \e\ is small enough then 
$$ \dim(\ker(\tilde{\go}^{n-2}:H^2(\tilde{X}) \into
H^{2n-2}(\tilde{X}))) = \dim(\ker(\go^{n-2}:H^2(X) \into
H^{2n-2}(X))) - 1$$
\end{lem}
\begin{proof}
By conveniently choosing $v_2$ and $v_0$ in \eqref{E:h2},
\begin{equation}\tag{\ref{E:h2}}
\begin{aligned}
(f^*(\go) + \e a)^{n-2} (f^*(v_2) + a v_0) =&\\
= \underset{H^{2n-2}(X)}{\underbrace{f^*(\go^{n-2}v_2 - \e^{n-2} t v_0)}}& + \e^{n-3}a^{n-2}((n-2)
\sigma v_0 + \e( \underset{\neq 0}{\underbrace{i^*v_2}} - c_1 v_0))
\end{aligned}
\end{equation}
the term in $H^{2n-2}(X)$ can be made equal to any
pre-chosen element in $\im(\go^{n-2}) \oplus \mbox{span}\{t\}$. Once
$v_2$ and $v_0$ are chosen, changing $v_2$ by an element in
$\ker(\go^{n-2})$ does not affect the result. On the other hand, by
varying $v_2$ by an element in $\ker(\go^{n-2})$ the coefficient of $a$ can be made
equal to anything in $H^{2}(M)$. Therefore $\dim(\im(\tilde{\go}^{n-2})) =
\dim(\im(\go^{n-2})) + 2$ and $\dim(H^2(\tilde{X})) = \dim(H^2(X)) +
1$, hence the result follows.
\end{proof}

Finally, we finish the study of the blow-up along surfaces claiming
that, for $i>2$, 
$$\dim(\ker(\tilde{\go}^{n-i})) = \dim(\ker({\go}^{n-i})).$$
Indeed, if $v_i \in \ker(\go^{n-i})$ then, $i^*(v_i) =0$, since it has
degree greater than 2, and therefore $a f^*(v_i) =a i^*(v_i) =0$ and 
$$(f^*(\go)+\e a)^{n-i}f^*(v_i)  = f^*(\go^{n-i} v_i) =0,$$
so $f^*(\ker(\go^{n-i})) \subset \ker(\tilde{\go}^{n-i})$.

Conversely, assuming $i$ even (the odd case is analogous),

\begin{align*}
(f^*(\go)+\e a)^{n-i}&(f^*(v_{i}) + a^{\frac{i}{2}} v_0 +
 a^{\frac{i-2}{2}}v_2) =\\
 &f^*(\go^{n-i}v_i) +
\e^{n-i-1} a^{n-\frac{i}{2}-1} ((n-i)\sigma v_0 + \e v_2) + \e^{n-i}
a^{n-\frac{i}{2}}v_0,
\end{align*}
and therefore $f^*(v_i) + a^{\frac{i}{2}} v_0 + a^{\frac{i-2}{2}}v_2$ will be in
$\ker(\tilde{\go}^{n-i})$ if, and only if, $v_0=0$ (by the coefficient
of $a^{n-\frac{i}{2}}$), $v_1 =0$ (by the coefficient of $a^{n-\frac{i}{2}-1}$) and $v_i
\in \ker(\go^{n-i})$, establishing the reverse inclusion.

So we have proved:
\begin{theo}{lefschetz in m2}
Let $i: M^2 \hookrightarrow X^{2n}$ be a symplectic embedding, $M$ and $X$ 
be compact 
and $\tilde{X}$ the blow up of $X$ along $M$. Then, for \e\ small enough,
\begin{itemize}
\item for $i >2$,
$$ \dim(\ker(\tilde{\go}^{n-i}))
= \dim(\ker(\go^{n-i}))
,$$
in particular, Lefschetz holds at level $i$ in $\tilde{X}$ if, and
only if, it does so in $X$;

\item if there is an element in $\ker(\go^{n-2})$ that restricts to a
nonzero element in $H^{2}(M)$ then
$$ \dim(\ker(\tilde{\go}^{n-2}))
= \dim(\ker(\go^{n-2}))
-1,$$
otherwise these kernels have the same dimension;

\item if there are elements $v_1,v_2 \in \ker(\go^{n-1})$ such that
$i^*(v_1 \wedge v_2) \neq 0$, then
$$ \dim(\ker(\tilde{\go}^{n-1}))
\leq \dim(\ker(\go^{n-1}))
-2,$$
otherwise
$$ \dim(\ker(\tilde{\go}^{n-1}))
 \leq \dim(\ker(\go^{n-1}))
.$$

\end{itemize}
\end{theo}

\subsection{The General case}
Now we treat the general case of the blow-up. Our main objective is to prove that if both $M^{2d}$ and $X^{2n}$
satisfy the Lefschetz property so does the blow-up of $X$ along $M$,
although in the course of this proof we obtain slightly  more, including a
generalization of Lemma \ref{T:h2}. The first part of the proof was already
encountered at the end of the 2-dimensional case.

\begin{prop}{i>2d}
Assume that $(M^{2d},\sigma) \hookrightarrow (X^{2n},\go)$ is a
symplectic embedding with $M$ and $X$ compact and $2d < n$. Let $\tilde{X}$ be the blown-up
manifold. Then, for $ i> 2d$
$$\dim(\ker(\tilde{\go}^{n-i})) = \dim(\ker(\go^{n-i})).$$
In particular, $\tilde{X}$ will satisfy the Lefschetz property at level
$i > 2d$ if, and only if, $X$ does so.
\end{prop}

\begin{remark}
The condition $2d < n$ is there only so that we can talk about a Lefschetz map at
level $i> 2d$, and this proposition says that {\em we can not change the
dimension of the kernel of the Lefschetz map beyond the dimension of the
submanifold along which we are blowing-up.}
\end{remark}

\begin{proof}
First we return to our usual notation and let $k =
n-d$. Let $v_i \in \ker(\go^{n-i}) \subset H^{i}(X)$ and consider the
cohomology class $f^*(v_i) \in H^i(\tilde{X})$. The restriction of $v_i$ to $M$ is zero,
since the degree of $v_i$ is greater than the dimension of
$M$. Therefore $a v_i= 0$ and
$$(f^*(\go) + \e a)^{n-i} f^*(v_i) = f^*(\go^{n-i} v_i) = 0.$$
Thus, $f^*(\ker(\go)^{n-i}) \subset \ker(\tilde{\go}^{n-i})$.

On the other hand assume that $v= f^*(v_i) + a v_{i-2} + \cdots + a^{l} v_{i-2l}$ 
is an element of the kernel of $\tilde{\go}^{n-i}$. We may
further assume that the last term above, $v_{i-2l}$, is not zero or else
$v$ is of the form $f^*(v_i)$. From $v \in \ker(\tilde{\go})$ we have
\begin{align*}
0 &= (f^*(\go) + \e a)^{n-i}(f^*(v_i) + a v_{i-2} + \cdots + a^{l} v_{i-2l})\\
&= f^*(\go^{n-i} v_i) + \sum_{j=0,m=1}^{n-i,l} \e^{j}
\begin{pmatrix}
n-i\\j
\end{pmatrix}
a^{j+m}
\sigma^{n-i-j} v_{i-2m}.
\end{align*} 
Since $i> 2(n-k)$, the degree of the element above is $ 2n -i < 2k$
and therefore the highest power of $a$ in the expression above is
still smaller than $k$. Hence the coefficient of $a^{l+n-i}$, which is 
$v_{i-2l}$, must vanish. Thus we had from the beginning $v=f^*(v_i)$
and the expression above reduces to
$$0 = (f^*(\go) + \e a)^{n-i}(f^*(v_i)) = f^*(\go^{n-i}v_i)$$
and $v \in f^*(\ker(\go^{n-i}))$, which shows the reverse inclusion
and proves the proposition.
\end{proof}

\begin{prop}{i=2d}
Assume that $i:(M^{2d},\sigma) \hookrightarrow (X^{2n},\go)$ is a
symplectic embedding with $M$ and $X$ compact and $2d < n$. Let $\tilde{X}$ be the blown-up
manifold. If there is a $v \in \ker(\go^{n-2d})$ such that $i^*v \neq
0$, then
$$\dim(\ker(\tilde{\go}^{n-2d})) =  \dim(\ker(\go^{n-2d})) -1,$$
otherwise these kernels have the same dimension, as long as \e\ is
small enough. In particular, if $X$ 
has the Lefschetz property at level $2d$, so does $\tilde{X}$.
\end{prop}
\begin{proof}
Initially we observe that the same argument used in Lemma
\ref{T:lemma2} shows that the existence of $v \in \ker(\go^{n-2d})$
such that $i^*v \neq 0$ is equivalent to the fact that the Thom class, $t$, 
of the embedding is not in the image of $\go^{n-2d}$. Now we let $k= n-d$
and  write down the Lefschetz map at level $2d$
\begin{equation}\label{E:i=2d}
\begin{aligned}
(f^*(&\go) + \e a)^{n-2d}(f^*(v_{2d})+a v_{2d-2} + \cdots + a^d v_0) = f^*(\go^{n-2d}v_{2d} - \e^{n-2d} v_0 t) +\\
&+ \sum_{i=k-d}^{k-1}
a^i\left(\left( \sum_{l \geq i - n+2d}^{d} \begin{pmatrix} n -2d\\i-l
\end{pmatrix} \e^{i-l}\sigma^{n-2d-i+l} v_{2(d-l)}\right)-  \e^{n-2d} v_0
c_{k-i}\right), 
\end{aligned}
\end{equation}
where the $c_i$'s are the Chern classes of the normal bundle of $M$. Then we
claim that we can make it equal to any element in 
\begin{equation}\tag{$*$}
f^*(\im(\go^{n-2d})\oplus\mbox{span}\{t\})\oplus a^{k-d} H^{2d}(M) \oplus \cdots
a^{k-1} H^{2}(M).
\end{equation}

The idea is the following: the system above is triangular
and therefore easy to solve.
Indeed, let $f^*(w_{2(n-d)})+a^{k-d} w_{2d} +\cdots + a^{k-1} w_{2}$ be an
element of the space ($*$). We start by choosing $v_{2d}$ and $v_0$ so that $\go^{n-2d}v_{2d}
- 
\e^{n-2d} t v_0$ equals $w_{2(n-d)}$. Observe that we can still change
$v_{2d}$ by any element in the kernel of $\go^{n-2d}$. Now look at the
coefficient of 
$a^{k-1}$ in \eqref{E:i=2d}:
$$ \e^{n-2d} v_2 + (n-2d) \e^{n-2d-1} \sigma v_0 - \e^{n-2d} v_0 c_1:=\e^{n-2d}v_2 + 
F(v_0).$$
Since we have already chosen $v_0$, we can now choose $v_2$ so that
the expression above equals $w_2$.

Assuming by induction that $v_{2j}$ have already been chosen for $j < j_0 < k-d$ so
that the coefficient of $a^{k-j}$ is $w_{2j}$ we see that the coefficient
of $a^{k-j_0}$ in \eqref{E:i=2d} is of the form
$$ \e^{n-2d} v_{2j_0} + F(v_0, \cdots, v_{2j_0-2}),$$
where $F$ is a function. Then again we can choose $v_{2j_0}$ so as to
have the desired equality.

Finally the coefficient of $a^{k-d}$ is of the form
$$ \e^{n-2d} i^*v_{2d} + F(v_0, \cdots, v_{2d-2}) \in H^{2d}(M).$$
And then, changing $v_{2d}$ by a multiple of the element in
$\ker(\go^{n-2d})$ whose restriction to $M$ is nonvanishing, we can make this
coefficient equal $w_{2d}$.

Now a simple counting of the dimensions involved shows that
$$\dim(\ker(\tilde{\go}^{n-2d})) = \dim(\ker(\go^{n-2d}))-1.$$

In order to prove the ``otherwise'' case, we start observing that if
for every $v\in \ker(\go^{n-2d})$, $i^*v =0$ then $f^*(\ker(\go)) \subset
\ker(\tilde{\go}^{n-2d})$. Therefore we immediately 
have $\dim(\ker(\go^{n-2d})) \leq \dim(\ker(\tilde{\go}^{n-2d}))$.

The reverse inequality is similar to what we have done so far and also to the subject of Proposition \ref{T:i<2d}, so we shall omit its proof.
\end{proof}

Before we can tackle the case $i < 2d$ we have to recall from Yan \cite{Ya96}.
\begin{lem}\label{T:yan}
 If $(M^{2d},\sigma)$ satisfies the Lefschetz property, there is a splitting of every cohomology class into primitive elements:
\begin{equation}\label{E:primitive splitting}
H^{i}(M) = P_i \oplus \im(\sigma),
\end{equation}
where $P_i$ is defined by
$$P_i = \{v \in H^{i}(M) \mid \sigma^{d-i+1}v =0 \},$$
if $i\leq d$ and $P_i = \{0\}$ otherwise. The elements in $P_i$
are called {\em primitive $i$--cohomology classes}.
\end{lem}

Hence, we can write every $v \in H^{i}$, $i\leq j$ in a unique way as $v = v^0 + v^1
\sigma + \cdots + v^{[i/2]}\sigma^{[i/2]}$, with $v^j$
primitive. Observe that if $i>d$, then the first few terms in this
decomposition will vanish simply because $P_j = \{0\}$ for
$j>d$. Again, the notation for the splitting above will be used
consistently in the sequence.

\begin{prop}{i<2d}
Let $i: (M^{2d},\sigma) \hookrightarrow (X^{2n},\go)$ be a symplectic
embedding with $M$ and $X$ compact and $2d < n$. Assume further that $M$ satisfies the
Lefschetz property. Then, for \e\ small enough and $i \leq 2d$,
$$\dim(\ker(\tilde{\go}^{n-i})) \leq \dim(\ker(\go^{n-i})).$$ 
In particular, if $X$ satisfies the Lefschetz property at level $i$ so does
$\tilde{X}$.
\end{prop}
\begin{proof}
Firstly we observe that the cases of $i$ odd and $i$ even
can be treated similarly, but for simplicity we shall work out only the 
even
case: $2i$. 

We want to take the limit $\e \into 0$ in the map
$\tilde{\go}^{n-2i}$, but, as it stands, the
resulting map will clearly have a big kernel. So, what we shall do is
to find linear maps $A_{\e}$ and $B_{\e}$ such that $\lim_{\e \into 0}
B_{\e}\tilde{\go}^{n-2i}A_{\e}$ has kernel $f^*(\ker(\go^{n-2i}))$. From this we shall
conclude that the dimension of the kernel of $\tilde{\go}^{n-2i}$ is
at most the dimension of the kernel of $\go^{n-2i}$ as long as \e\ is
small enough.

We define $A_{\e}:H^{2i}(\tilde{X}) \into H^{2i}(\tilde{X})$ by
$$A_{\e}\left(f^*(v_{2i}) + \sum_{j=0}^{i-1} a^{i-j} v_{2j}\right) =
f^*(v_{2i})  + 
\sum_{j=0}^{i-1} \frac{1}{\e^j}a^{i-j} v_{2j}.$$
And $B_{\e}: H^{2n-2i}(\tilde{X}) \into H^{2n-2i}(\tilde{X})$ by
$$B_{\e} \left(f^*(v_{2n-2i}) + \sum_{j=0}^{i-1} a^{n-d-i+j}
v_{2d-2j}\right) =f^*(v_{2n-2i}) + \sum_{j=0}^{i-1}
\frac{1}{\e^{n-d-2i+j}}a^{n-d-i+j} v_{2d-2j} $$

Now we move on to write the map $\lim_{\e\into 0}B_{\e} \tilde{\go}^{n-2i}
A_{\e}$:
$$\lim_{\e\into 0}B_{\e} \tilde{\go}^{n-2i}
A_{\e} = f^*(\go^{n-2i}v_{2i}) + \sum_{j=0}^{i-1}
a^{n-d-i+j}\sum_{l=0}^{i-1}b_{d-j-l}\sigma^{d-j-l}v_{2l},$$
where $b_j= \begin{pmatrix} n-2i \\ j \end{pmatrix}$ are the binomial
coefficients.

We can further split the cohomology classes $v_{2l}$ into their primitive
parts, according to lemma \ref{T:yan}, $v_{2l} = v_{2l}^0 + \sigma
v_{2l}^1 + \cdots +\sigma^l v_{2l}^l$. With that, elements of
$H^{2i}(\tilde{X})$ will be in the kernel of the map above only if the
coefficients of $a^j\sigma^l$ vanish. The only terms that will give us information about primitives of degree $2l$ are the coefficients of $a^{k-i+l} \sigma^{d-2l}, a^{k-i+l+1} \sigma^{d-2l-1}, \dots, a^{k-1}\sigma^{d-l-i+1}$, and the vanishing of these is equivalent to the following:
$$\begin{pmatrix}
b_{d-2l} & b_{d-2l-1} & \cdots & b_{d-l-i+2} &  b_{d-l-i+1}\\
b_{d-2l-1} & b_{d-2l-2} & \cdots & b_{d-l-i+1} &  b_{d-l-i}\\
\vdots&       & \ddots & &\vdots\\
b_{d-2l-i+2} & b_{d-2l-i+1} & \cdots & b_{d-2i+4} &  b_{d-2i+3}\\
b_{d-2l-i+1} & b_{d-2l-i} & \cdots & b_{d-2i+3} &  b_{d-2i+2}
\end{pmatrix}
\begin{pmatrix}
v_{2l}^0\\
v_{2l+2}^1\\
\vdots\\
v_{2i-4,}^{i-2-l}\\
v_{2i-2}^{i-1-l}
\end{pmatrix}
=0$$
in the case $2i <d$, and a similar matrix for $2i>d$. What is important here is that in both cases the matrix will be constant along its anti-diagonals (it is a Toeplitz matrix) and the top right entry is nonzero.
Now, if we can prove that all the matrices above are invertible, we
will conclude that $f^*(v_{2i}) + \sum a^{i-j} v_{2j}$ is in the
kernel of $\lim B_{\e}\tilde{\go}^{n-2i}A_{\e}$  if and only if
$v_{2j}=0$ for all $j < i$ and $v_{2i} \in
\ker\{\go^{n-2i}\}$. So the next lemma finishes the theorem.
\noqed\end{proof}

\begin{lem}\label{T:oliver}
Let $b^n_j = \begin{pmatrix} n\\ j \end{pmatrix}, n, j \in \N$. Then for any $p \in \N$
$$\Delta^{n,p+1}_k = \det \begin{pmatrix}
b^n_{k+p} & b^n_{k+p-1} & \cdots & b^n_{k+1} &  b^n_{k}\\
b^n_{k+p-1} & b^n_{k+p-2} & \cdots & b^n_{k} &  b^n_{k-1}\\
\vdots&       & \ddots & &\vdots\\
b^n_{k+1} & b^n_{k} & \cdots & b^n_{k-p+2} &  b^n_{k-p+1}\\
b^n_{k} & b^n_{k-1} & \cdots & b^n_{k-p+1} &  b^n_{k-p}
\end{pmatrix}
\neq 0
$$
if $b^n_k \neq 0$. 
\end{lem}
\begin{proof}
Initially we observe that $b^n_k \neq 0$ if and only if $n\geq k \geq 0$ and for $n=k$ the matrix above has zeros above the anti-diagonal and ones on it, so the determinant is a power of $-1$. Further, by adding to each column the one to its right and using the binomial identity $b^n_k + b^n_{k-1} = b^{n+1}_k$ we get
$$\Delta^{n,p+1}_k = \det \begin{pmatrix}
b^{n+p}_{k+p} & b^{n+p-1}_{k+p-1} & \cdots & b^{n+1}_{k+1} &  b^n_{k}\\
b^{n+p}_{k+p-1} & b^{n+p-1}_{k+p-2} & \cdots & b^{n+1}_{k} &  b^n_{k-1}\\
\vdots&       & \ddots & &\vdots\\
b^{n+p}_{k+1} & b^{n+p-1}_{k} & \cdots & b^{n+1}_{k-p+2} &  b^{n}_{k-p+1}\\
b^{n+p}_{k} & b^{n+p-1}_{k-1} & \cdots & b^{n+1}_{k-p+1} &  b^n_{k-p}
\end{pmatrix}
$$
Now it is easy to check that 
$$\Delta^{n+1,p+1}_k  =\frac{(n+p+1)!(n-k)!}{n!(n+p-k+1)!}\Delta^{n,p+1}_k,$$
showing that $\Delta^{n+1,p+1}_k$ is nonzero if $\Delta^{n,p+1}_k$ is nonzero and we obtain the result by induction.
\end{proof}

These three propositions give us the following

\begin{theo}{lefschetz}
Let $i: (M^{2d},\sigma) \hookrightarrow (X^{2n},\go)$ be a symplectic
embedding with $M$ and $X$ compact and both satisfying the
Lefschetz property and $2d < n$. Let $(\tilde{X},\go+\e \ga)$ be the blow-up of $X$ along $M$ with the symplectic form from Theorem \ref{T:cohomology}. Then, for \e\ small enough, $\tilde{X}$ also satisfies the Lefschetz property.
\end{theo}

\section{Massey Products and the Blow-up}\label{S:massey}

Having determined how the Lefschetz property behaves under blow-up, we turn our attention to formality.
Here, we use Massey products to prove that manifolds are not formal, since formality implies that these products vanish.  The object this section is to prove that under mild codimension conditions, Massey products are preserved in the blow-up. This will allow us to find examples of nonformal symplectic manifolds next section.

The ingredients for a Massey product are  $a_{12},a_{23},a_{34} \in
\gO^{\bullet}(M)$, three closed forms such that $a_{12}a_{23}$ and $a_{23}a_{34}$ are exact. Then, denoting $\bar{a}= (-1)^{|a|}a$, we define
\begin{equation}\tag{$*$}
\begin{cases}
\overline{a_{12}}a_{23} &= da_{13}\\
\overline{a_{23}}a_{34} &= da_{24}.
\end{cases}
\end{equation}
In this case, one can consider the element
$\overline{a_{13}} a_{34} + \overline{a_{12}} a_{24}$. By the choice of $a_{13}$ and $a_{24}$ this form is closed, hence it represents a cohomology class. Observe,
however, that $a_{13}$ and $a_{24}$ are not well defined and we can change
them by any closed element, hence the expression above does not define a
unique cohomology class but instead an element in the
quotient $H^{\bullet}(M)/\mathcal{I}([a_1],[a_3])$, where $\mc{I}$ denotes the ideal generated by its arguments.

\begin{definition}
The {\it triple Massey product} or just {\it triple product}
$\langle [a_{12}], [a_{23}], [a_{34}]\rangle$, 
of the cohomology classes $[a_{12}]$, $[a_{23}]$ and $[a_{34}]$
with $[a_{12}] [a_{23}]=[a_{23}] [a_{34}]=0$
 is the coset
$$\langle [a_{12}],[a_{23}],[a_{34}]\rangle = [\overline{a_{12}}
  a_{24} + \overline{a_{13}} a_{34}] + ([a_{12}], [a_{34}])
  \in H^{\bullet}(\mc{A})/\mathcal{I}([a_{12}],[a_{34}]),$$
where $a_{13}$ and $a_{24}$ are defined by $(*)$. 
\end{definition}

\begin{theo}{massey products}
Let $i:M^{2(n-k)} \hookrightarrow X^{2n}$ be a symplectic embedding  with
$M$ compact
and let $\tilde{X}$ be the blown-up manifold, then:
\begin{itemize}
\item if $X$ has a nontrivial triple Massey product, so does $\tilde{X}$,
\item {\em (Babenko and Taimanov \cite{BT00a})} if $M$ has a nontrivial triple Massey product and $k > 3$, so does $\tilde{X}$.
\end{itemize}
\end{theo}
\begin{proof}
We start with the first claim and assume the Massey product $\langle v_1, v_2, v_3
\rangle$ is nonzero
in $X$. This means that there is $u$ representing such a product
with $[u] \not \in \mathcal{I}([v_1], [v_3])$, the ideal generated by $[v_1]$ and
$[v_3]$ in $H^*(X)$. If we consider the
product $\langle f^*v_1, f^*v_2,f^*v_3 \rangle$ we see that $f^*u$ is a
representative for it. The question then is whether $f^*[u]$ is in
the ideal $(f^*[v_1], f^*[v_3])$. Let us assume there was a relation
of the kind 
$$f^*[u] =  f^*[v_1](f^*\xi_1+a \zeta_1^1 + \cdots a^{k-1}
\zeta_1^{k-1}) +  f^*[v_3](f^*\xi_3+a \zeta_3^1 + \cdots a^{k-1}
\zeta_3^{k-1})$$  
Then, using the product rules \eqref{E:product rules},
$$f^*[u] = f^*([v_1]\xi_1 +[v_3] \xi_3) + a (i^*[v_1] \zeta_1^1 +  i^*[v_3]
\zeta_3^1)+ \cdots + a^{k-1}  (i^*[v_1] \zeta_1^{k-1} + i^*[v_3] \zeta_3^{k-1}).$$
Now, since the sum above is a direct one, all the
coefficients of the powers of $a$ must vanish and the following must hold:
$$f^*[u] = f^*( [v_1] \xi_1 + [v_3] \xi_3).$$ 
Since $f^*$ is an injection, we conclude that $[u] \in ([v_1], [v_3])$
which contradicts our initial assumption. 

Now we treat the second case. We
start by assuming that $v_1$, $v_2$ and $v_3 \in \gO(M)$ are closed forms
satisfying
$$ v_1 \wedge v_2 = dw_1 ~~~\mbox{ and }~~~ v_2\wedge v_3 = dw_2,$$
with $[w_1 v_3 -(-1)^{|v_1|} v_1 w_2] \not \in ([v_1],[v_3])$.
Letting $\gf:\tilde{V} \into V$ be the map of diagram \eqref{E:diagram}
and $\pi:V \into M$ the projection of the disc bundle, we have the
following relations in $H^*(\tilde{X})$
$$ \ga \gf^*\pi^* v_1 \wedge \ga \gf^*\pi^* v_2 = d (\ga^2
\gf^*\pi^*w_1) ~~~ \mbox{ and } ~~~  \ga \gf^*\pi^* v_2 \wedge \ga
\gf^*\pi^* v_3 = d (\ga^2 \gf^*\pi^*w_2)$$
The question then is again whether the cohomology class of the form
$$\ga \gf^*\pi^*v_1 \ga^2 \gf^*\pi^* w_1-(-1)^{|v_1|}\ga \gf^*
\pi^*v_1 \ga^2 \gf^*\pi^*w_2$$
is in the ideal generated by $a [v_1]$ and $a [v_3]$.

Suppose it was. Then there would be a relation of the type
\begin{align*}
a^3  [w_1 v_3 - (-1)^{|v_1|}v_1 w_2] &= a[v_1](f^*\xi_1+a \zeta_1^1+
\cdots + a^{k-1} \zeta_1^{k-1}) +\\& + a[v_3](f^*\xi_3+a \zeta_3^1+
\cdots + a^{k-1} \zeta_3^{k-1})\\
&= a([v_1] i^*\xi_1 [v_3] i^*\xi_3)+ a^2([v_1] \zeta_1^1 + [v_3]
\zeta_3^1) + \cdots +\\ &+ a^{k-1} ([v_1] \zeta_1^{k-2} + [v_3]
\zeta_1^{k-2}) + a^k ([v_1] \zeta_1^{k-1} + [v_3] \zeta_3^{k-1}).
\end{align*}
Expanding $a^k$ and using again that the result is a direct sum, we
look at the coefficient of $a^3$. Comparing both sides we see that it equals $[w_1 v_3 -
(-1)^{|v_1|}v_1 w_2]$, so
$$[w_1 v_3 - (-1)^{|v_1|}v_1 w_2] = [v_1] \zeta_1^2 + [v_3]
\zeta_3^2 - c_{k-3} \zeta_1^{k-1} [v_1] + c_{k-3} \zeta_3^{k-1} [v_3].$$
But this contradicts the fact that  $[w_1 v_3 -(-1)^{|v_1|} v_1 w_2]
\not\in ([v_1],[v_3])$.
\end{proof}

\section{Examples}\label{S:nonformal examples}

In this section we give concrete examples where the blow-up procedure can be used to create manifolds satisfying the Lefschetz property. The examples we produce will also have nontrivial Massey products, therefore producing a counter-example to the conjecture of Babenko and Taimanov. 

Our starting point are nilmanifolds, i.e., compact quotients of a nilpotent Lie group by a maximal lattice. It is a result of Benson and Gordon \cite{BG88} that nontoroidal nilmanifolds never satisfy the Lefschetz property. Also, Nomizu's theorem \cite{No54} implies that the Lie algebra of the correspoding Lie group with its differential $(\wedge^{\bullet}\Gg^*,d)$ furnishes a minimal model for the nilmanifold, therefore, no nontoroidal nilmanifold is formal. Indeed, it is a result of Cordero {\it et al} \cite{CFG86} that they always have nontrivial (maybe higher order) Massey products.

The simplest nilmanifold with the properties we need is the one obtained from the product of two copies of the Heisenberg group.

\begin{ex}{heisenberg}
If $G$ is the 3 dimensional Heisenberg group then the Lie algebra has a basis 
formed by the left invariant vector fields whose values at the
identity are 
$$
\del_1 = \begin{pmatrix}
0&1&0\\
0&0&0\\
0&0&0\\
\end{pmatrix};~~~
\del_2 = \begin{pmatrix}
0&0&0\\
0&0&1\\
0&0&0\\
\end{pmatrix};~~~
\del_3 = \begin{pmatrix}
0&0&-1\\
0&0&0\\
0&0&0\\
\end{pmatrix}
$$

Then we check that $[\del_1,\del_2] = -\del_3$ and $[\del_1, \del_3] = [\del_2,\del_3] =0$.
Therefore, the quotient manifold, $\mathbb{H}$, of the 3--Heisenberg group
by the lattice generated by $\exp \del_i$ will have Chevalley-Eilenberg 
complex $(\wedge^{\bullet}\frak{g}^*,d)$  generated
by the invariant 1-forms $e_i$ dual to the $\del_i$ related by
$$ de_1 = de_2 =0;~~~de_3 = e_1\wedge e_2.$$
So, this manifold is the nilmanifold associated to the Lie algebra with structure $(0,0,12)$.

Hence $H^1(\mathbb{H})$ is generated by $\{e_1,e_2\}$,
$H^2$, by $\{e_{13}, e_{23}\}$ and $H^3$, by $\{ e_{123}\}$, where, as usual, $e_{ij}$ is the shorthand for $e_i\wedge e_j$. Therefore $b_1= b_2= 2$ and $b_0 = b_3 =1$
\vskip6pt
\noindent 
{\bf Massey products.}
With the relations above for $e_1$, $e_2$ and $e_3$ we get that
$$ e_1 \wedge e_2 = de_3 ~~~~ \mbox{ and } ~~~~ e_2 \wedge e_1 =
d(-e_3).$$

Therefore we can form the Massey product $\langle e_1, e_2, e_1 \rangle = 
e_3 \wedge e_1
+ e_1 \wedge (-e_3) = - 2 e_{13} \neq 0$. Observe that in this
case, since $e_1 \wedge e_2$ is exact, Massey products have no
indeterminacy and the above is a nontrivial one in $\mathbb{H}$.
\end{ex}

\begin{ex}{gil}
Now consider the product $\mathbb{H} \times \mathbb{H} = (0,0,12,0,0,45)$.
The triple product $\langle e_1, e_2,e_1 \rangle$ is still nonzero. Further, the form
$$\go = e_{14} + e_{23} + e_{56}$$
is closed and has top power $6 e_{123456}$, which is
everywhere nonvanishing. Hence $\go$ is a symplectic form in
$\mathbb{H} \times \mathbb{H}$.

It is easy to see
that the kernel of $\go:H^2 \into H^4$ is $\mbox{span}\{e_{25}\}$ and the kernel of $\go^2$ in $H^1$ is $\mbox{span}\{e_2,e_5\}$.

In $\mathbb{H}$, consider the path
$$\ga(t) = \exp(t(\del_1+\del_2+\del_3)) =
\begin{pmatrix}
1 & t & t+ \dfrac{1}{2} t^2\\
0 & 1 & t \\
0 & 0 & 1
\end{pmatrix}, ~~t \in [0,2],$$
then $\ga' = \del_1 + \del_2 + \del_3$ and, besides this, $\ga(2) \approx \ga(0)$,
hence this is a circle.

In $\mathbb{H} \times \mathbb{H}$ there are two copies of $\ga$ (one in
each factor) making 
a torus $T^2$. A basis for the tangent space of this torus is given by
$\{\del_1+\del_2+\del_3, \del_4 + \del_5 + \del_6\}$. The symplectic form
evaluated on this basis equals 1 everywhere, hence this torus is a
symplectic submanifold. On the other hand, $e_{25}$
evaluated on this basis also equals 1 everywhere.

Therefore, by Theorems \ref{T:lefschetz in m2} and \ref{T:massey products}, the blow-up $M^6$ of $\mathbb{H} \times \mathbb{H}$ 
along this torus satisfies the Lefschetz property (for \e\ small enough)
and has a nontrivial triple product.
\end{ex}

\begin{ex}{munoz}
Still let $M$ be the manifold from the previous example. We can change the symplectic form slightly so that it is rational and still satisfies the Lefschetz property. Hence an appropriate multiple of it will represent an integral cohomology class and hence, by Donaldson's theorem \cite{Do96}, it is Poincar\'e dual to a symplectic submanifold $(N^4,\omega)$. Using Fern\'andez and Mu\~noz' result on formality of Donaldson submanifolds \cite{FM02}, $N$ still satisfies the Lefschetz property and by Donaldson's theorem the inclusion $N \hookrightarrow M$ induces an isomorphism $H^1(M) \cong H^1(N)$ and an injection $H^2(M) \hookrightarrow H^2(N)$. Now, the Massey product in $M$ comes from three 1-forms and therefore still exists in $N$ and further, since we have an injection in $H^2$, this product is nonzero in $N$. So $N$ is a nonformal symplectic 4-manifold satisfying the Lefschetz property.
\end{ex}

\begin{ex}{gil2}
Let $(N^4,\sigma)$ be the manifold obtained in Example \ref{T:munoz}
that has a nontrivial triple product and satisfies the Lefschetz property.
By construction, $\sigma$ is an integral cohomology class, therefore, $(N, \sigma)$ can be symplectically embedded in
$\mathbb{C}P^{6}$, by Gromov's Embedding Theorem \cite{Gr71,Ti77}. By Theorem
\ref{T:lefschetz}, the blow-up of $\mathbb{C}P^{6}$ along $N$ will have the
Lefschetz Property. According to Theorem \ref{T:massey products}, it will have a
nonvanishing triple product (and thus is not formal) and from Theorem
\ref{T:cohomology} it is 
simply connected.
\end{ex}

\bibliographystyle{abbrv}
\bibliography{references}

\begin{thebibliography}{10}

\bibitem{BT00b}
I.~K. Babenko and I.~A. Taimanov.
\newblock {\em Massey products in symplectic manifolds}.
\newblock Sb. Math, 191:1107--1146, 2000.

\bibitem{BT00a}
I.~K. Babenko and I.~A. Taimanov.
\newblock {\em On nonformal simply connected symplectic manifolds}.
\newblock Siberian Math. J., 41:204--217, 2000.

\bibitem{BG88}
C.~Benson and C.~S. Gordon.
\newblock {\em K\"ahler and Symplectic Structures on Nilmanifolds}.
\newblock {Topology}, 27:513--518, 1988.

\bibitem{Br88}
J.~Brylinski.
\newblock {\em A differential complex for symplectic manifolds}.
\newblock {J. Differential Geometry}, 28:93--114, 1988.

\bibitem{Ca04b}
G.~R. Cavalcanti.
\newblock {\em New aspects of the $dd^c$-lemma}.
\newblock PhD thesis, Oxford University, 2004.

\bibitem{CFG86}
L.~Cordero, M.~Fern\'andez, and A.~Gray.
\newblock {\em Symplectic manifolds with no {K}\"ahler structure}.
\newblock {Topology}, 25:375 -- 380, 1986.

\bibitem{DGMS75}
P.~Deligne, P.~Griffiths, J.~Morgan, and D.~Sullivan.
\newblock {\em Real homotopy theory of {K}\"ahler manifolds}.
\newblock {Invent. Math.}, 29:245--274, 1975.

\bibitem{Do96}
S.~Donaldson.
\newblock {\em Symplectic submanifolds and almost-complex geometry}.
\newblock {J. Differential Geometry}, 44:666--705, 1996.

\bibitem{FM02}
M.~Fern\'andez and V.~Mu{\~n}oz.
\newblock {\em Formality of {D}onaldson submanifolds}.
\newblock math.SG/0211017, to appear in Math. Zeit.

\bibitem{Go95}
R.~E. Gompf.
\newblock {\em A new construction of symplectic manifolds}.
\newblock {Ann. of Math. (2)}, 142:527--595, 1995.

\bibitem{Gr71}
M.~Gromov.
\newblock {\em A topological technique for the construction of 
  differential equations and inequalities}.
\newblock In {\em Actes Congr. Internat. Math.}, volume~2, pages 221--225, Nice
  1970, 1971.

\bibitem{Gu03}
M.~Gualtieri.
\newblock {\em Generalized Complex Geometry}.
\newblock PhD thesis, Oxford University, 2003.
\newblock math.DG/0401221.

\bibitem{Hi03}
N.~Hitchin.
\newblock {\em Generalized {C}alabi-{Y}au Manifolds}.
\newblock {Quart. J. Math. Oxford}, 54:281--308, 2003.

\bibitem{IRTU03}
R.~Ib{\'a}{\~n}ez, Y.~Rudyak, A.~Tralle, and L.~Ugarte.
\newblock {\em On certain geometric and homotopy properties of closed
  symplectic manifolds}.
\newblock In {\em Proceedings of the Pacific Institute for the Mathematical
  Sciences Workshop ``Invariants of Three-Manifolds'' (Calgary, AB, 1999)},
  volume 127, pages 33--45, 2003.

\bibitem{Kos85}
J.~L. Koszul.
\newblock {\em Crochet de {S}chouten-{N}ijenhuis et cohomologie}.
\newblock In {\em Elie Cartan et les meth. d'aujourd'hui}, pages 251 -- 271.
  Ast\'erisque hors-s\'erie, 1985.

\bibitem{LO94}
G.~Lupton and J.~Oprea.
\newblock {\em Symplectic manifolds and formality}.
\newblock {J. Pure Appl. Algebra}, 91:193 -- 207, 1994.

\bibitem{Ma90}
O.~Mathieu.
\newblock {\em Harmonic cohomology classes of symplectic manifolds}.
\newblock {Comment. Math. Helv.}, 70:1--9, 1990.

\bibitem{McD84}
D.~McDuff.
\newblock {\em Examples of simply-connected symplectic non-{K}\"ahlerian
  manifolds}.
\newblock {J. Differential Geometry}, 20:267--277, 1984.

\bibitem{MS00}
D.~McDuff and D.~Salamon.
\newblock {\em Introduction to Symplectic Topology}.
\newblock Oxford Mathematical Monographs. Oxford University Press, 1995.

\bibitem{Me98}
S.~A. Merkulov.
\newblock {\em Formality of canonical symplectic complexes and {F}robenius
  manifolds}.
\newblock {Internat. Math. Res. Notices}, 14:727--733, 1998.

\bibitem{Mi79}
T.~J. Miller.
\newblock {\em On the formality of $k-1$ connected compact manifolds of
  dimension less than or equal to $4k-2$}.
\newblock {Illinois J. Math.}, 23:253--258, 1979.

\bibitem{No54}
K.~Nomizu.
\newblock {\em On the cohomology of compact homogeneous spaces of nilpotent Lie
  groups.}
\newblock {Ann. of Math. (2)}, 59:531--538, 1954.

\bibitem{Pa66}
A.~N. Parshin.
\newblock {\em A generalization of the {J}acobian variety {\em (Russ.)}}.
\newblock Isvestia, 30:175--182, 1966.

\bibitem{RT00}
Y.~Rudyak and A.~Tralle.
\newblock {\em On {T}hom spaces, {M}assey products and non-formal symplectic
  manifolds}.
\newblock {Internat. Math. Res. Notices}, 10:495--513, 2000.

\bibitem{Ti77}
D.~Tischler.
\newblock {\em Closed 2--forms and an embedding theorem for 
  manifolds}.
\newblock {J. Differential Geometry}, 12:229--235, 1977.

\bibitem{Ya96}
D.~Yan.
\newblock {\em Hodge Structure on Symplectic Manifolds}.
\newblock {Adv. in Math.}, 120:143--154, 1996.

\end{thebibliography}

\end{document}